# From the theory of "congeneric surd equations" to "Segre's bicomplex numbers"

CINZIA CERRONI[1]

**Abstract:** We will study the historical pathway of the emergence of *Tessarines* or *Bicomplex numbers*, from their origin as "imaginary" solutions of irrational equations, to their insertion in the context of study of the algebras of hypercomplex numbers.

## 1. Introduction

Beginning from the end of the first half of the nineteenth century, particularly in Great Britain, in the wake of researches on geometrical interpretation of complex numbers, studies developed that led to the birth of new systems of hypercomplex numbers and are at the basis of the birth of modern algebra. In particular, the discovery in 1843 of Quaternions by William Rowan Hamilton (1805-1865) revealed to mathematicians the existence of an algebraic system that had all the properties of real and complex numbers except commutativity of multiplication. Further, the studies on symbolic algebra[2] by George Peacock (1791-1858) and on logic by Augustus De Morgan (1806-1871) created a context of reflection and analysis on the laws of arithmetic and their meaning.

As a result, researches were carried out in Great Britain on new systems of hypercomplex numbers, leading to the discovery of Octonions (1843/1845) by John T. Graves (1806-1870) and Arthur Cayley (1821-1895), of the theory of Pluriquaternions[3] (1848) by Thomas Penyngton Kirkman (1806-1895), of Biquaternions (1873) and the Algebras (1878) of William Kingdon Clifford (1845-1879) and to the systematic presentation of the algebras of hypercomplex numbers (1870), known until that time, through the work of Benjamin Peirce (1809-1880).

---

[1] Università degli Studi di Palermo, Dipartimento di Matematica e Informatica. E-mail: cinzia.cerroni@unipa.it
[2] For an analysis of the rise of symbolic algebra in Great Britain see [Pycior 1981], [Pycior 1982].
[3] [Kirkman 1848].

In this framework, *tessarines* (or bicomplex numbers) were identified. The idea came to James Cockle beginning from the observation made by William G. Horner on the existence of irrational equations, called *"congeneric surd equations",* which admit neither real solutions nor complex solutions. Specifically, beginning from the equation $\sqrt{j}+1 = 0$, Cockle introduced a new imaginary unity $j$, which satisfies $j^2 = 1$, and taking inspiration from Hamilton's theory of quaternions wrote the generic tessarine as $w = a + bi + cj + ijd$, where $j^2 = 1$. He then immediately hypothesized a generalization of it, introducing *octrines* at the end of his work and conjecturing generalization to sixteen unities. In a series of subsequent articles he analyzed the algebraic properties of tessarines, noticing the existence of zero divisors, and compared these properties with those of other systems of hypercomplex numbers known at the time, like *Quaternions*, *Octonions* etc., also in reference to the symbolic algebra introduced by Peacock. In 1892 Corrado Segre, in his complex geometry studies, and in the general context of the study of the algebras of hypercomplex numbers set going by Karl Weierstrass (1815-1897) in 1884, rediscovered the algebra of bicomplex numbers. Segre presented bicomplex numbers as the analytical representation of the points of bicomplex geometry and recognized that Hamilton himself in the study of biquaternions had introduced the same quantities; but neither Hamilton before, nor Segre after, noticed Cockle's works on the subject. Segre studied the algebra of bicomplex numbers, also determining a decomposition of them, which is a particular case of Peirce's.

In this paper, we will analyze the historical pathway of the emergence of tessarines or bicomplex numbers, from their origin as "imaginary" solutions of irrational equations, to their insertion in the context of study of the algebras of hypercomplex numbers.

## 2. "Congeneric surd equations"

The theory of irrational equations, referred to as "congenerics", was worked out by William George Horner (1786-1837)[4] in 1836 and came into being in the context of the study of irrational equations, when it was observed that irrational equations exist that have neither real solutions nor complex solutions. The considerations on this subject that W. G. Horner wrote in a letter to Thomas

---

[4]He was a British mathematician. He was educated at the Kingswood School, Bristol, where he became an assistant master at the age of fourteen. After four years he was promoted to headmaster. In 1809 Horner left Bristol to found his own school at 27 Grosvenor Place in Bath. Horner's only significant contribution to mathematics lies in the method of solving algebraic equations which still bears his name. Contained in a paper submitted to the Royal Society (read by Davies Gilbert on 1 July 1819), "A new method of solving numerical equations of all order by continuous approximation", it was published in *Philosophical Transactions of the Royal Society of London* (Ser. A, 1 (1819), pp.308-335). This method was ascribed to him by A. de Morgan and J. R. Young. Actually Paolo Ruffini (1765-1822) anticipated Horner's method in 1802 (Mem. Coronata Soc. Ital. Sci., 9 (1802), pp. 44-526). See http://www-history.mcs.st-and.ac.uk/Biographies/Horner.html.

Stephens Davies (1794?-1851),[5] were published by the latter. Specifically, in the introduction to the work Davies wrote:

*"If those mathematicians who have met with a quadratic equation[6] whose 'roots' either under a real or imaginary form could not be exhibited will recall to memory the surprise with which they viewed the circumstance, and the attempts which they made to solve the mystery, they will read with no ordinary gratification the following discussion of the general question of which this forms a part. The general theory of such equations, very happily named by Mr. Horner 'Congeneric Equations' is here laid down with great clearness, and, so far as I know, for the first time,- as it is, indeed, nearly the first time the formation of any general and philosophic views respecting them has been attempted [...]."* [T. S. Davies in Horner 1836, p. 43].

From the considerations by T. S. Davies it emerges that the need was felt to go deeper into the matter, in that it had not yet received the necessary attention and the analysis undertaken by W. G. Horner was finally seen by Davies as clarifying it.

In this connection, in his article W. G. Horner, after some considerations on the need to study the case in which irrational equations do not admit any real root, wrote as follows:

*"My Dear Sir, I agree with you in thinking that the properties of irrational equations have not received that degree or kind of attention from writers on the elements of algebra, which was due either to the importance of the subject, or to a consideration for the comfort of young students. […] In solving equations involving radicals every one has experienced the necessity of putting his results to the proof before he could venture to decide which of them, or whether any of them, could be trusted; but as the latter alternative, or the failure of **every** result, is of rare occurrence in books of exercises […] In the mean time even classical writers have spoken of clearing an equation from radicals, in order to its solution, as a process of course, and which would not in any way affect the conditions. The consequence is, that a habit prevails of talking about equations without any regard to this peculiar case, and therefore in language which when applied to it becomes quite incorrect."*

[Horner 1836, pp. 43-44]

---

[5] He was a British mathematician. He became a fellow of the Royal Society of Edinburgh in 1831, and he contributed several original and elaborate papers to its *Transactions*. He also published "Researches on Terrestrial Magnetism" in the *Philosophical Transactions*, "Determination of the Law of Resistance to a Projectile" in the *Mechanics' Magazine*, and other papers in the *Cambridge and Dublin Mathematical Journal*, the *Civil Engineer*, the *Athenæum*, the *Westminster Review*, and *Notes and Queries*. In April, 1833 he was elected a Fellow of the Royal Society. In 1834, he was appointed one of the mathematical masters at the Royal Military Academy at Woolwich. His death, after six years of illness, took place at Broomhall Cottage, Shooter's Hill, Kent, on 6 January 1851. See https://en.wikipedia.org/wiki/Thomas_Stephens_Davies.

[6] He refers to quadratic irrational equations, which he exemplified as follows: "*For instance,* $2x + \sqrt{x^2 - 7} = 5$, *the 'roots' of which are* 4 *and* $\frac{8}{3}$ *are determined by the common process; neither of which substituted in the equation reduces it to zero. These are the roots of its congeneric surd equation* $2x - \sqrt{x^2 - 7} = 5$.

He analyzed the reasons why it was not correct to apply, without conditions, the "standard" procedure of simplification of radicals: *"For, being aware that if **unequal** quantities (+ a, - a) be raised to power, denoted by equal exponents, the power **may** nevertheless be equal; we are assured that, conversely, if of equal quantities roots be extracted which are denoted by equal exponents, the roots may nevertheless be **unequal.** This remark furnishes a sufficient reason for rejecting the third pair of principles and consequently the ordinary method of clearing an equation from surds."* [Horner 1836, p. 46] and proposed that the possible roots of the given irrational equation should be sought among the solutions of a "group of equations" constituted by the given irrational equation, by the associated "congeneric" irrational equation and by the rational equation obtained by simplifying the radicals of the given equation irrational by raising to power or alternatively as a product of the irrational equation and its "congeneric." Specifically, he affirmed: *"Thus, between the statements $a = \sqrt{x}$ or $a - \sqrt{x} = 0$ and $a^2 = x$ or $a^2 - x = 0$ the copula $a + \sqrt{x} = A$ has been lost sight of. The complete chain is*

$$a - \sqrt{x} = 0 \quad \text{or} \quad A$$
$$a + \sqrt{x} = A_1 \quad \text{or} \quad 0$$
$$a^2 - x = 0.$$

*You are well aware that this copula will, in all cases of surd equations, consist of all the variations that can be made of the given formula by varying the affection of each radical it contains in all possible ways. […] (for the sake of a convenient term I would venture to say, its congeners). […] We know that the continued product of a surd formula and all its congeners will produce a rational formula; and that such rational formula, being equated to zero, may be solved by as many roots as it has dimensions. We are also certain that each of these roots will cause one of the congeneric surd formulae to vanish; otherwise the product of all would not be = 0 as assumed."* [Horner 1836, p. 49]

Horner, therefore, wondered what the meaning was of finding a "solution" to the rational equation associated with the irrational equation, and if this solution could be called the "root" of the given irrational equation or not: *"But, is the value of x which effects this, to be called a root of the surd formula? No, it is a root of the rational combination only. – Have irrational equations, then, no roots?"* [Horner 1836, p. 49] and arrived at the conclusion that in the case of "congeneric" irrational equations, instead of looking for the "roots", it was necessary to look for the "solutions" to the group of equations constituted by the irrational equation, by its congeneric and by the associated rational equation: *"None at all. – What have they, then, in the place of roots? An equitable chance, in common with each formula in the congeneric society, of solution by means of the solution of the stock-equation."* [Horner 1836, p. 49].

The article concludes with some questions posed by the author, regarding how to set problems on this type of irrational equations and on what degree to attribute to the latter: "*But if an equation has no root, nor even a certainty of solutions, in what form can it be intelligibly proposed? A note of interrogation posed either for solution or correction. – To what order can surd equations be assigned? To the fractional order $\frac{m}{n}$, when n congeneric formulae produce a rational equation of the mth order? […] Are the chances of solution equal for each individual congener?*" [Horner 1836, pp. 49-50] .

The importance of Horner's article does not lie in the observations made – though these are interesting, seeing that nobody had treated the matter in any depth until that moment – but rather in having inspired J. Cockle's works on "Tessarines."

## 3. The "Tessarines" and the systems of quadruple Algebra of James Cockle[7]

Between 1848 and 1850 James Cockle (1819-1895), in a series of works introduced and studied a new system of hypercomplex numbers, which he himself called "*Tessarines.*" The idea arose from W. Horner's considerations on "*congeneric*" irrational equations. In this connection, Cockle wrote as follows to T. S. Davies: "*You are aware that the new algebraical symbol, which I have recently discovered, was suggested to my mind by reflecting on the structure of congeneric surd equations. Did I require any excuse for throwing the present investigations into a form of a letter to you, I should find it in the fact that the remarks of Horner on congenerics were put into a similar shape.*"[Cockle 1848, p. 435]

Davies himself, in the introduction to Cockle's work, which he himself proposed for publication, underlined the central role played by Horner's considerations: "*I believe that Mr. Horner's attention was first directed by me to the consideration of what he has happily termed 'congeneric surd equations'. This was done in consequence of a difficulty that I had found in discussing such an equation which was accidentally brought before me; a difficulty which, from the friendship that had long subsisted between us, I very naturally referred to him. His reply was printed in the Philosophical Magazine about a dozen years ago; and it has received proper attention from those algebraists who look with due care to the fundamental principles of their*

---

[7] James Cockle was born in 1819 in Great Oakley, Essex, England and died in London in 1895. He studied at Trinity College, Cambridge. He published papers of both pure and applied mathematics and of history of science. He was involved in fluid dynamics, magnetism, algebra and differential equations. He was a member of the Royal Astronomical Society, of the Cambridge Philosophical Society, of the London Mathematical Society and of the Royal Society in London. See http://www-history.mcs.st-and.ac.uk/Biographies/Cockle.html.

*science. Almost every new difficulty in algebra is the precursor of more extended views respecting its principles; and often out of a case which may, at first sight, appear to be little more than a conundrum, very important improvements in our general theories often arise. All effort to generalize our views, so as to include the 'conundrum', deserve, therefore, attentive consideration; and for this reason I send you for publication the accompanying letter which I have just received from your correspondent and my friend, Mr. James Cockle, of the Middle Temple.* […]" [Davies in Cockle 1848, p. 435].

He observed that a new difficulty in algebra often proves to be the precursor of an extension of its principles, as, precisely, in the case of Horner's "congeneric equations", which inspired the introduction of "*Tessarines.*"

J. Cockle, within the paradigm of the symbolic algebra introduced by George Peacock and stimulated by the discovery of Hamilton's *Quaternions*, introduced a new "symbol" that he called "impossible quantity", as the solution of an "impossible" irrational equation, and beginning from it introduced the *Tessarines*.

Indeed, Cockle having been trained mathematically in the school of Peacock, wrote: "[…] *but to express my opinion that Dr. Peacock's Algebra is the standard work on the philosophy of that science, and that his views form a near approximation to those which will ultimately be received as forming the true basis of symbolical algebra. His view appears to form a basis somewhat analogous to the geometrical one which MONGE has taken in his principle of "contingent relations".* [Cockle 1848a, p. 365]", he identified in algebraic symbolism the paradigm of the discipline and hence identified in the need to introduce "impossible quantities", the result of the laws of algebra: "[…] *yet an algebraist of the SCHOOL OF PEACOCK sees, in the occurrence of those quantities, an inevitable result of the laws imposed upon his symbols, and, consequently, a fit subject for his researches.*" [Cockle 1848a, p. 365].

Specifically, he proceeded as follows. Inspired by Horner's works and therefore by reflections on the fact that some irrational equations exist which can be defined "impossible" since they admit neither real nor imaginary solutions, he believed it was necessary to introduce some quantities that "would represent" such equations: "*I think that I am justified in stating that there is a species of algebra which has as yet received no attention from the cultivators of analysis; I mean the algebra of impossible quantities, or impossible expressions; for the term quantity is not strictly applicable to such cases. In speaking of impossible quantities, I must be understood as alluding to something very different from the unreal or imaginary quantities of ordinary algebra.* […] *Premising, then, that there are algebraic equations the very supposition of the existence of which involves an arithmetical contradiction – equations which have no root whatever – I shall designate* *such equations as impossible. Bearing this designation in mind: By an impossible quantity is meant a root of an impossible equation."* [Cockle 1848a, p. 364-365].

Hence an "impossible quantity" is a root of an impossible equation. The introduction of a "symbol" that describes the "simplest" element of the "impossible quantity", in analogy with imaginary unity, therefore represented a development of algebra originating from the paradigm of symbolic algebra: *"[...] But the impossible quantities which we are about to consider, are inevitably forced upon our notice in the contemplation of expressions which ordinary algebra presents to us. Hence the admission of a symbol to denote impossible quantity – should such symbol ever be admitted – into works on algebra, would be an extension of and perhaps an innovation upon the previously existing science, and would, in fact, be a legitimate development of it."*

[J. Cockle 1848a, p. 365]

Hence Cockle introduced the new "symbol" as the solution of the simplest "impossible" irrational equation, in analogy with positive and negative numbers and imaginary unity; specifically:

*"Let us now proceed to give some kind of system to this algebra of impossibles. The ordinary algebra takes cognizance of three species of quantity; - positive, negative, and unreal. Perhaps the simplest forms under which these three kinds of quantity can be exhibited are the following*:

$$+1 \;;\; -1 \;;\; (-1)^{1/2} \;;$$

*which we may call*

$$p \;;\; n \;;\; \text{and } u \;;$$

*respectively. Now what is the corresponding simple representation of impossible quantities?*

*The process which I shall employ for the purpose of answering this question is as follows. I shall express p, n, and u by equations, and then, guided by the analogy afforded by them, I shall form a corresponding equation for an impossible quantity which I shall denote by i, and which will probably be the simplest equation by means of which impossible quantity can be exhibited.*

*Thus, p is the root of the equation*

$$0 = 1 - x \;;$$

*so, n is the root of*

$$0 = 1 + x \;;$$

*also u is a root of*

$$0 = 1 + x^2 \;;$$

*and guided by analogy, I take i to be the root of*

$$0 = 1 + \sqrt{x} \;.$$ [Cockle 1848a, p. 366].

Subsequently, after introducing the simplest "impossible quantity" as the solution of the simplest "congeneric" irrational equation, influenced by Hamilton's works on *Quaternions*, he introduced *Tessarines*, as follows:

"*Let*

$$1 + i^2 = 0$$

*then i is the simplest representation of unreal quantity; and, if*

$$1 + \sqrt{j} = 0,$$

*then j is the simplest representation of impossible quantity. It is to borne in mind that, in the latter equation, the radical is to be considered as essentially affected with the sign +.*

*Let w, x, y, z be any real quantities, positive, negative, or zero; also let*

$$w + ix + jy + ijz = t$$

*then t I call a tessarine, and w, x, y, z its constituents. The latter term I have adopted from the quaternion theory of Sir W. R. Hamilton*". [Cockle 1848, p. 436].

Again inspired by Hamilton's theory, he introduced the "equations" that describe the *Tessarines*, that is to say the relationships that intervene between imaginary unity and the "impossible quantity", as follows:

"*For convenience, I shall denote the product of i and j by k; I shall assume that $j^2$ is equal to unity. I say assume because, although I have reasons (not, however, free from objection), for making such a supposition, yet I wish to reserve myself full liberty to modify the assumption, and to discuss ab initio the symbol $j^n$ and its symbolic relation to j. Then the following system of relations will, together with the second equation given in this letter, furnish us with all the conditions requisite for the formation of a theory. That system is*

$$i^2 = k^2 = -1 \quad j^2 = 1$$
$$ij = k, \, jk = i, \, ki = -j.\text{"}$$

[Cockle 1848, p. 437]

In introducing the relationships between the new symbol *j* and imaginary unity *i*, Cockle moves away from considerations on "impossible quantities", through which he had introduced it, and in actual fact generalizes the procedure, to the point that the article ends up hypothesizing the possibility of introducing "*Octrines*" and "beyond": "*By way of conclusion, I will add that, if we had three independent imaginaries i, j, k, we should have to deal, not with tessarines, but with what might be called octrines. These last mentioned expressions would be of the form*

$$w + ix + jy + kz + ijp + ikq + jkr + ijks,$$

*where w, x, y, z, p, q, r, s are any real quantities. These quantities might be termed the constituents of the octrine, which could not vanish unless all its constituents were zero. So, if we had four independent imaginaries [...] the corresponding expression would consist of sixteen terms, and would have properties analogous to those of the tessarine and the octrine. I must here break off [...]"* [Cockle 1848, pp. 438-439].

From the latter reference, it emerges that the author had jauntily generalized the procedure, probably also influenced by the considerations on Octonions by John Graves and Arthur Cayley[8].

After introducing the *Tessarines*, Cockle studied their properties: "[…] *The next step is, to ascertain the fundamental properties of the new symbol, its origin and nature being duly considered.*" [Cockle 1849, p. 38]

Specifically he showed that a *Tessarine* is annulled if its constituents are equal to zero [Cockle 1848, p. 436], that the product of two *Tessarines* is a *Tessarine* [Cockle 1848, p. 437], [Cockle 1849a, p. 408], and that, as we say today, there are zero divisors. Cockle defined the latter property "anomalous": "*But I am about to point out another anomaly,* […]: *it is that on the supposition that $j^2$ is equal to the unity*,

$$(1 - j) \times (1 + j) = 1 - 1 = 0;$$

*this is to say, the zero may be considered as the product of two impossible factors, neither of which vanishes.* […] *It may be said, is zero, then decomposable into non-vanishing factors? Impossible. I reply, true, the factors are impossible: they are so by their origin and nature.*" [Cockle 1849, p. 42], and he justified it himself, precisely because of the impossible nature of *Tessarines*. It is to be observed that this was the first time that the possibility of the existence of zero divisors was stressed. This is a property that at first, as observed, appeared anomalous to the author, but subsequently was to be approved by him, faced with the fact that *tessarines* satisfy the other laws of ordinary algebra.

Indeed, Cockle continued to explore "the nature" of *Tessarines* comparing them with Hamilton's *Quaternions*, and identifying the properties that from his point of view made them adhere to the laws of "ordinary" algebra. Specifically he observed that while in the case of *Quaternions* the characteristic relationships force us to sacrifice the commutative law, in the case of *Tessarines* this is not the case: "*My present object is, however, connected with Mr. Boole's observation on the 'Laws of Quaternions'* […] *He has there shown that the quaternion relations*

$$i^2 = j^2 = k^2 = -1$$
$$ij = k, jk = i, ki = j,$$

*when considered as of universal application conduct us to the conditions*

---

[8] For in-depth examination of the birth and development of Octonions cf. [Baez, 2001], [Cerroni, Vaccaro 2010], [Cerroni 2010], [Gray, Parshall 2007].

$$ji = -k,\ kj = -i,\ ik = -j,$$

and thus lead to a sacrifice of the commutative character of multiplication. It becomes therefore a subject of inquiry whether the fundamental equations of my Tessarine System entail upon us the necessity of any such sacrifice. To ascertain this, I shall apply the method of Mr. Boole to the Tessarine conditions

$$i'^2 = -j'^2 = k'^2 = -1$$
$$i'j' = k',\ j'k' = i',\ k'i' = -j'$$

Let, then, the subject of operation be $j'y$, and we have

$$i'j'j'y = k'j'y$$

or
$$i'j'^2 y = k'j'y$$

but
$$j'^2 = 1,$$

therefore
$$i'y = k'j'y.$$

So, if the subject be $k'z$, we have

$$j'k'k'z = i'k'z,$$

or
$$j'k'^2 z = -j'z = i'k'z;$$

and, proceeding thus, we see that the commutative character of multiplication is preserved in the Tessarine Theory." [Cockle 1849f, p. 124]

From this analysis it emerges that Cockle had realised that the commutative property of multiplication in *Tessarines* is a consequence of the characteristic relationships between imaginary unities and therefore that it is the relationships between imaginary unities that determine the "nature" of systems of hypercomplex numbers. Hence Cockle continued in the analysis of possible relationships between imaginary unities and introduced two new algebras, *Coquaternions* and *Cotessarines*: "If α, β, and γ be three imaginaries of which the respective squares are equal either to positive or negative unity, and which are subject to the relation $αβ = γ$, then these imaginaries afford us four and only four, essentially distinct systems of Quadruple Algebra. The first system in that of QUATERNIONS; the second, that of TESSARINES; for the third system I propose the name of COQUATERNIONS; and for the fourth that of COTESSARINES." [Cockle 1849g, p. 197]

He proceeded as follows. For each case he considered the squares of the imaginary unities and the relationship $αβ = γ$ and multiplying both members of the latter by the imaginary unities obtained the characteristic relationships of the four algebras. In this way he recovered the *Quaternions* and the *Tessarines* and identified two new algebras, *Coquaternions* and *Cotessarines*.

Specifically, he examined the case in which two of the three imaginary unities have squares equal to unity, obtaining the system of the *Co-quaternions*:

"3. The Third, or Coquaternion System.

*In this case, the characteristic equations may be expressed by*

$$-a^2 = b^2 = c^2 = 1 \ldots\ldots\ldots(15.),$$

where a, b, and c are the imaginaries. We have, also, from $\alpha\beta = \gamma$,

$$ab = c \ldots\ldots\ldots\ldots(16.)$$

*To avoid repetition, I shall now use the terms multiplier and multiplicand with reference to the numbered equation next preceding those terms wherever they are used.*

*Let a be the multiplier, we have*

$$aa\,b = a^2 b = ac,$$

*or by (15.),*

$$-b = ac \ldots\ldots\ldots\ldots(17.)$$

*Take c as multiplicand [...] we find*

$$-bc = a \ldots\ldots\ldots\ldots(18.)$$

*With b as multiplier, we find*

$$-c = ba \ldots\ldots\ldots\ldots(19.);$$

*and with a as multiplicand, we obtain*

$$ca = b \ldots\ldots\ldots\ldots(20.)$$

*So, c as multiplier gives,*

$$a = cb \ldots\ldots\ldots\ldots(21.)$$

*Hence, in this Coquaternion System, the commutative character of multiplication is lost, and the system is, consequently, abnormal. [...]"* [Cockle 1849g, p. 198].

And finally, considering the case in which the three imaginary unities have a unitary square, he obtained the system of the *Cotessarines*:

*"4. The Fourth, or Cotessarine System.*

*Let d, e, and f be the imaginaries of this system, then its characteristic is that*

$$d^2 = e^2 = f^2 = 1,$$

*and we also take de = f. Proceeding as suggested [...] we find, successively,*

$$e = df, \text{ and } d = fe;$$

$$ef = d, \text{ and } fd = e;$$

*and, finally, f = ed. This system is normal. One of its defects is, that it takes no cognizance of the expression $\sqrt{-1}$, and consequently has no relation to ordinary Double Algebra."*

[Cockle 1849g, p. 198].

Comparing the four algebras Cockle established which of them are a *normal* system[9], i.e. those that satisfy the laws of ordinary algebra, and which are a *non-normal* system. Specifically,

---

[9] This nomenclature was introduced by Cockle himself

*Tessarines*, unlike *Quaternions*, are *normal* and for this reason, according to the author, they are the most adequate to extending the ordinary algebra of couples: "[…] *we see that the Tessarine System is normal.* […] *This Tessarine System appears to me to be the natural extension of ordinary Double Algebra.*" [Cockle 1849g, p. 198].

In a series of works[10] in which he studied the structure of *Tessarines*, Cockle determined the form of their module and its properties. Subsequently he compared the module of Tessarines with that of Cotessarines, Quaternions and Coquaternions: "[…] *where I have shown the existence of these four systems, discussed them, and pointed out their characteristic. I have proposed to advert for a moment to the same subject, to consider it under a slightly different aspect, and also to exhibit, for convenient of comparison, the modular expressions of all the systems.*" [Cockle 1849b, pp. 434-435]. Other works in which there are further observations are [Cockle 1850a], [Cockle 1850b]. Cockle's last work on these matters dates from 1852 [Cockle 1852].

From the previous considerations, it emerges that Cockle and the English context in which he was immersed had a "modern" conception of the laws of algebra and their nature and consequently jauntily introduced "new" systems of hypercomplex numbers, motivating them both intrinsically with the laws of algebra and with their possible applications to geometry.

## 4. W. R. Hamilton's biquaternions

As long ago as 1850, William Rowan Hamilton (1805-1865) had developed an extension of quaternions, defining the algebra of biquaternions. In that year, he made it the object of a communication to a meeting of the British Association for the Advancement of Science in Edinburgh. All that is extant is the following Report:

"*the author briefly explained the term which he had been obliged to introduce into this new system; showed the simplicity and the reasons for the leading operations in it; and by a few very simple experiments on the rotation of planes round axes inclined to each other, explained the simple interpretation of some of those results which appeared at first to be inconsistent with the principles of the ordinary analysis.*" [Hamilton, 1852]

In 1853 he wrote "Lectures on Quaternions" [Hamilton 1853], in which he collected all his results on the theory of Quaternions and extended them. In the introduction, he noticed that his researches had inspired many contemporaries, showing that he was aware of developments and connections with his theory, in the context of studies in Great Britain. Specifically, he affirmed:

---

[10] [Cockle 1849c], [Cockle 1849d], [Cockle 1849e].

"*My thanks are due, at this last stage, to the friends who have cheered me throughout by their continued sympathy; to the scientific contemporaries* [note: *In these countries, Messrs Boole, Carmichael, Cayley, Cockle, De Morgan, Donkin, Charles and John Graves, Kirkman, O'Brien, Spottiswoode, Young, and perhaps others: some of whose researches or remarks on subjects connected with quaternions* (*such as the triplets, tessarines, octaves, and pluriquaternions*) *have been elsewhere alluded to, but of which I must regret the impossibility of giving here a fuller account.*[…]] *who have at moments turned aside from their own original researches, to notice, and in some instance to extend, results or speculations of mine*; […]" [Hamilton 1853, Preface p. 64].

And subsequently, after mentioning some studies by contemporaries in § XCIV of the index of the contents, he writes: "[…] *subsequent extension (in the same year) by J. T. Graves, Esq., to a theorem respecting* **sums of eight squares**, *and to a theory of certain* **octaves**, *involving seven distinct imaginaries*; *allusion to subsequent publications of Professor De Morgan, and other mathematicians of these countries, in the same general field of research, or at least on analogous subjects, such as* **triplets**, **tessarines**, *and* **pluriquaternions**; […]" [Hamilton, 1853 p. liii]. In the paragraphs mentioned, he observed that it was not possible for him to go deeper into all those results: "[…] *But it is impossible for me here to attempt to do any kind of justice to the talent and candour of the many able and original mathematical writers in these countries, who have been pleased to acknowledge that some subsequently published speculations of theirs, on subjects having some general connexion with or affinity to the present one, were, more or less, suggested or influenced by the quaternions.*[…]" [Hamilton, 1853 p. 539].

Hence it emerges that Hamilton knew and appreciated the studies by his contemporaries and particularly Cockle's studies on "*Tessarines.*" In going deeper into the algebra of quaternions, and in particular the theory of equations with coefficients in quaternions, in his "Lectures on Quaternions" Hamilton again introduced "*biquaternions*, as imaginary solutions of quadratic equations in quaternions:

"*We see however, that the imaginary solutions of the proposed equation in quaternions still present themselves under the* GENERAL FORM,

$$q = q' + \sqrt{-1}\ q''$$

*where q' e q'' are real quaternions, while* $\sqrt{-1}$ *is still the old and the ordinary imaginary of algebra, and is distinguished from all those other roots of negative unity which are peculiar to the present calculus,* […] *and II$^{st}$, by its being, as a factor, commutative with every other. An expression of this general form is called by me* BIQUATERNIONS.*"* [Hamilton 1853, pp. 638-639].

He stressed the importance of considering the imaginary solutions of quadratic equations with coefficients in quaternions, for completeness of the theory, in analogy with the imaginary solutions of equations with real coefficients. Specifically, he affirmed: "*The theory of such **biquaternions** is as necessary and important a complement to the theory of **single** or **real** quaternions, as in algebra the theory of **couples**, or of expressions of the form*

$$x' + \sqrt{-1}\, x'',$$

*where x' and x'' denote some **two** positive or negative or null numbers, is to the theory of **single** or **real numbers** or quantities. It is admitted that the doctrine of **algebraic equations** would be entirely incomplete, if their **imaginary roots** […] were to be neglected, or kept out of view. And in like manner we may already clearly see, from the foregoing remarks and examples, that no theory of **equations in quaternions** can be considered as complete, which refuses or neglects to tale into account the **biquaternions solutions** that may exist, of the form above assigned, in any particular or general inquiry.[…]*"." [Hamilton 1853, p. 639].

It can hence be observed that Hamilton, in introducing biquaternions, followed a line of thought analogous to that of Cockle in the case of Tessarines. That is to say, new "hypercomplex numbers" are introduced, to complete the theory of equations with coefficients in quaternions.

He also gave some numerical examples in support of the theory. He considered the quadratic equation with coefficients in quaternions

$$q^2 = qi + j$$

and determined its six solutions, two of which are quaternions while four are biquaternions. Specifically, they are: $q = \frac{1}{2}(i-k) \pm \frac{1}{2}(i+j)$ and $q = \frac{1}{2}i\left(1 \mp \sqrt{-3}\right) - k$; $q = \frac{1}{2}(i+k) \pm \frac{1}{2}(1-j)\sqrt{-3}$, where $\sqrt{-3}$ is the imaginary unity [Hamilton 1853, p. 641].

Hamilton, in studying biquaternions, noticed that they admit zero divisors: "*It must, however, be confessed that such calculations as these with biquaternions, […] are sometimes very delicate, and require great caution, from the following circumstance […]. This circumstance is that the product of two biquaternions may vanish, without either factor separately vanishing. To give a very simple example, the product*

$$\left(k + \sqrt{-1}\right)\left(k - \sqrt{-1}\right) = k^2 + 1 = 0$$

*while $k + \sqrt{-1}$ and $k - \sqrt{-1}$ much each be considered as different from zero […].*" [Hamilton 1853, p. 650], and called biquaternions satisfying such properties "Nullific" or "Nullifiers": "*It seems convenient, therefore, to call biquaternions of this class Nullific or to say that they are Nullifiers […].*" [Hamilton 1853, p. 672].

In 1866, in *Elements of Quaternions*, a work published posthumously, Hamilton again returned to the subject, studying quaternions that he denominated "planar", that is to say complanar with one of the imaginary units, for instance *i* [Hamilton 1866, p.113; p. 240]. He observed that they can be represented in the form $x + iy$, with $x$ and $y$ real, therefore constituting a real sub-algebra isomorphic to the field of complex numbers [Hamilton 1866, p. 244]. Subsequently, he introduced "*complanar biquaternions*", as solutions of equations with coefficients in complanar quaternions, and showed that they can be put in the form $x_1 + hy_1 + i(x_2 + hy_2)$, with $x_1$, $y_1$, $x_2$ and $y_2$ real and *h* and *i* imaginary units [Hamilton 1866, p.277], that is in the form of those that C. Segre called "*bicomplex*" numbers, which are in fact Cockle's "*Tessarines*." Although he knew Cockle's works, Hamilton did not notice the analogy between "*complanar biquaternions*" and "*Tessarines*", an analogy noticed by Segre himself, as we will see in the next section. Hence, in Hamilton's formulation, bicomplex numbers constitute a commutative sub-algebra of the algebra of biquaternions. A few years later, in 1873[11] it was William Kingdon Clifford that introduced two other types of "biquaternions", again in the form $a + b\omega$, with *a* and *b* quaternions, but respectively with unity $\omega$ such that $\omega^2 = 0$ or $\omega^2 = +1$.[12] Other developments came in 1878, when Clifford, through a reflection on the comparison between Hamilton's quaternions and the vectorial language of Hermann Günther Grassmann, generalized this theory to those that are known as "*Algebras of Clifford*."[13] In this paper we will not deal all these works, except as regards the development of studies on bicomplex numbers.[14]

# 5. C. Segre's bicomplex numbers

In 1892, continuing a line of thought already begun some years before[15] [Segre 1889-90, 1890], Corrado Segre (1863-1924) published the work "*Real representations of complex forms and hyperalgebraic bodies*" [Segre 1892] in which he inserted the geometrical interpretation of the algebra of bicomplex numbers, returning after forty years to the interrupted thread of Hamilton's thought: "[…] *following that principle of the extension of notions, to which Mathematics owes so much progress, and that in particular for the geometry of algebraic varieties had led from real points to complex points. Now a further extension appears appropriate. Complex points are not*

---

[11] [Clifford 1873].
[12] A comparison was made between Hamilton's "*biquaternions*" and Clifford's by Arthur Buchheim in 1885 [Buchheim 1885, p. 293].
[13] [Clifford 1878].
[14] For an in-depth examination of the history of geometric calculation see [Freguglia 2004].
[15] On this aspect of Segre's work see [Brigaglia 2013], [Zappulla 2009].

*enough. It is useful to introduce some bicomplex points, that is to say some entities that have through images the complex points of representative forms.[16]"* [Segre 1892, p. 449].

After dealing at length with the hyperalgebraic entities (complex entities that in a real representation are algebraic), Segre introduced bicomplex points as a natural completion of the complex projective straight line. In the last part of the work he introduced and studied, precisely, a new kind of numbers, which are the analytical representation of bicomplex points, bicomplex numbers: *"The introduction of imaginary points in geometry corresponds to the introduction of imaginary numbers (coordinates) in analysis. What will the further generalization be of the concept of number that will correspond to the extension that we have made of the geometric field by introducing bicomplex points?"*[17] [Segre 1892, p. 455].

He proceeded as follows. Let us suppose the complex points of a straight line are represented by the real points of the plane σ, and precisely the point that on the straight line has as its coordinate the complex number x+iy (x real y and $i^2$=-1) has as its image in the plane σ the point of coordinates (x,y). Then to obtain the bicomplex points on the same straight line, one must also consider in the plane σ the complex points (x,y), whose coordinate are $x = x_1 + hx_2$ and $y = y_1 + hy_2$, where $x_1$, $x_2$, $y_1$ and $y_2$ are real and $h^2 = -1$, and therefore consider that the bicomplex point of the straight line has the coordinate x+iy, that is to say in the plane

$$x_1 + hx_2 + i(y_1 + hy_2) = x_1 + hx_2 + iy_1 + hiy_2.$$

In this way on the straight line the bicomplex point will have as its coordinates *"bicomplex numbers"* of the type $x_1 + hx_2 + iy_1 + hiy_2$. where $x_1$, $x_2$, $y_1$, $y_2$ are real and *i* and *h are* two distinct imaginary units for which $h^2 = i^2 = -1$ (but $h \neq \pm i$) and their product is associative and commutative. The algebra of bicomplex numbers is therefore presented as the commutative algebra of the numbers *x + iy*, where *x* and *y* are complex numbers (in the imaginary unity, distinguished from *i*, denoted as *h*). Each bicomplex number is therefore expressed as a linear combination of the units 1, *i, h* and $k = ih$ with $i^2 = h^2 = -1$ and $k^2 = 1$. [Segre 1892, p. 456].

As mentioned in the previous section, it was Segre himself who realized that bicomplex numbers are equivalent to Hamilton's *"planar biquaternions"*: "[…] *but from the need to solve algebraic equations in which the unknowns are quaternions, being conducted then to consider alongside the quaternions, ordinary or real, whose 4 coefficients (scalar) are real, the imaginary quaternions or biquaternions in which those coefficients are imaginary, one obtains for the*

---

[16] *"[…] seguendo quel principio dell'ampliamento delle nozioni, a cui la Matematica deve tanti progressi, e che in particolare per la geometria delle varietà algebriche aveva portato dai punti reali ai punti complessi. Ora si presenta opportuna un'ulteriore estensione. Non sono più sufficienti i punti complessi. Conviene introdurre dei punti bicomplessi, cioè degli enti che abbiano per immagini i punti complessi delle forme rappresentative."*

[17] *"L'introduzione dei punti imaginari in geometria corrisponde all'introduzione dei numeri imaginari (coordinate) in analisi. Quale sarà l'ulteriore generalizzazione del concetto di numero, che corrisponderà all'estensione che abbiam fatta del campo geometrico introducendo i punti bicomplessi?"*

*biquaternions of a given plane the representation (Elements, no. 257) $x + iy = x_1 + h\, x_2 + i(y_1 + h\, y_2)$, where h is an imaginary number such that $h^2 = -1$. Apart from the distinction between the meanings of versor and number that Hamilton attributes respectively to the two symbols i and h, his biquaternions of a given plane are the same thing as our bicomplex numbers.*"[18] [Segre 1892, p. 457]. But Segre too, like Hamilton, did not realize that *bicomplex numbers* coincide with the *Tessarines* introduced by Cockle in 1848. Segre, moreover, also makes reference [Segre 1892] to Hermann Hankel (1839 – 1873), who in 1867[19] had hinted at the subject, referring to Hamilton's biquaternions, and to Rudolf Lipschitz (1832 – 1903), who in 1886[20] dealt more at length with the subject, referring to bicomplex numbers.

Segre, in particular, inserted the algebra of bicomplex numbers in the general framework of the study of hypercomplex algebras, which can be said to have been given systematic form with the publication, a few years earlier, in 1884, of the letter of Karl Weierstrass to Hermann Schwarz,[21] followed by the works of Schwarz himself[22] and above all of Richard Dedekind.[23] In the eighteen-eighties the subject became of particular interest in German mathematics and, above all thanks to Eduard Study, became an integral part of researches related to the groups of Lie.[24]

Segre's study on this algebra is ample and exhaustive. The first point highlighted is, naturally, the presence of zero divisors, which Segre called "*Nullifics*[25]" taking up the term used by Hamilton [Segre 1892 p. 458]. He determined that the zero divisors in bicomplex numbers constitute two ideals (in the modern sense), which he called "*infinite sets of nullifics*[26]", respectively generated by $h + i$ ($I_1$) and by $-h + i$ ($I_2$). Specifically he wrote: "*Those of the $1^{st}$ set are those bicomplex numbers that would be annulled if in place of the symbol i we had -h: they are the products of any bicomplex numbers [...] for the number h + i [...]. Those of the $2^{nd}$ set are those bicomplex numbers that would be annulled for i = h; that is to say, the products of any bicomplex numbers [...] for the number -h + I [...].*"[27] [Segre 1892 p. 459]. He also observed that the product

---

[18] "[…] *ma dai bisogni della risoluzione delle equazioni algebriche in cui le incognite sono quaternioni, essendo poi condotto a considerare accanto ai quaternioni, ordinari o reali, i cui 4 coefficienti (scalari) son reali, i quaternioni immaginari o biquaternioni in cui quei coefficienti sono immaginari, ottiene per i biquaternioni di un dato piano la rappresentazione (Elements, n. 257) $x + iy = x_1 + h\, x_2 + i(y_1 + h\, y_2)$, ove h è un numero imaginario tale che $h^2 = -1$. A parte la distinzione fra i significati di versore e di numero che Hamilton attribuisce risp. ai due simboli i ed h, i suoi biquaternioni di un dato piano son la stessa cosa che i nostri numeri bicomplessi.*"
[19] [Hankel 1867].
[20] [Lipschitz 1886 pp.125-131].
[21] [Weistrass 1884].
[22] [Schwarz 1884].
[23] [Dedekind 1885]. For the history of the theory of algebras in the eighteen-eighties see [Lutzen 2001].
[24] For in-depth examination in this direction see [Hawkins 2000].
[25] "*Nullifici*"
[26] "*schiere infinite di nullifici*"
[27] "*Quelli della 1° schiera sono quei numeri bicomplessi che si annullerebbero ove in luogo del simbolo i vi si ponesse –h: essi sono i prodotti di numeri bicomplessi qualunque [...] pel numero h + i [...]. Quelli della 2° schiera sono quei*

of two non-zero bicomplex numbers will be equal to zero if and only if they respectively belong to the *"two sets of nullifics"*[28], that is to say to the two ideals.

Segre determined the geometrical interpretation of the *"Nullifics"* in the projective geometry of bicomplex numbers. Specifically, he defined among the bicomplex varieties the two sets of *proto-strings*[29], (*"infinite strings"*[30]) on a (complex) straight line and showed that a *proto-string* will have as the coordinates of the points the bicomplex numbers of a coset module $I_1$ or $I_2$ depending on the type of proto-string involved: *"the points of the objective straight line that are on a same proto-string of the 1$^{st}$ set are those that have for coordinates different bicomplex numbers for nullifics of the 2$^{nd}$ set."*[31] [Segre 1892 p. 459].

Subsequently, he proceeded to determine the structure of bicomplex algebra. Specifically he determined two decompositions of it. He observed, in particular, that given a bicomplex number $x + iy$, setting $Z = x + hy$, $Z' = x - hy$ and $g = \dfrac{1-hi}{2}$, $g' = \dfrac{1+hi}{2}$, one obtains $x + iy = Zg + Z'g'$, with $g \in I_1$ and $g' \in I_2$. He therefore established that a bicomplex number can be decomposed into the sum of two zero divisors, each of the two possible types: *"every bicomplex number x + i y she can be decomposed in a clearly determined way into the sum of two nullifics, one* (Zg) *of the 1$^{st}$ set, the other* (Z'g') *of the 2$^{nd}$"*[32] [Segre 1892, p. 459].

Having established this first result, the Piedmontese mathematician was able to deduce another, more significant, decomposition. He proceeded as follows; given the bicomplex number $x + iy = x_1 + hx_2 + iy_1 + hiy_2$, setting $X = x_1 - y_2$, $Y = x_2 + y_1$, $X' = x_1 + y_2$, $Y' = -x_2 + y_1$, he obtained $x + iy = Xg + Yk + X'g' + Y'k'$, con *X, Y, X'* and *Y'* evidently real and $k = hg = \dfrac{h+i}{2}$, $k' = -hg' = \dfrac{-h+i}{2}$. From this he deduced, given $g^2 = g$, $k^2 = -g$, $kg = gk = k$ (and likewise for $g'$ and $k$),[33] in modern terms, that the bicomplex numbers ("nullifics") of the form $Xg + Yks$ constitute a sub-algebra isomorphic with ordinary complex numbers and therefore the algebra of bicomplex numbers is isomorphic with the direct sum of two copies of the field of complex numbers: *"The nullifics of the 1$^{st}$ set are therefore reduced to complex numbers* (with real coefficients X, Y) *with*

---

*numeri bicomplessi che s'annullerebbero per i = h; ossia i prodotti di numeri bicomplessi qualunque [...] pel numero – h + i [...]."*
[28] *"due schiere"*
[29] *protofili*
[30] *"infiniti fili"*
[31] *"i punti della retta oggettiva che stanno su uno stesso protofilo della 1° schiera sono quelli che hanno per coordinate numeri bicomplessi differenti per nullifici della 2° schiera"*
[32] *"ogni numero bicomplesso x + i y si può scomporre in un modo ben determinato nella somma di due nullifici l'uno* (Zg) *della 1° schiera, l' altro* (Z'g') *della 2°"*
[33] That is to say, *g* and *k* (*g'* and k') satisfy the same relationships as 1 and *i*.

*the two unities g, k; and likewise the nullifics of the 2<sup>nd</sup> set with the two unities g', k' […]"*[34] [Segre 1892 p. 460]. Segre immediately determined the geometrical meaning of the real coefficients $X$, $Y$ and $X'$, $Y'$ of the new representation of the bicomplex number, that is to say: *"the new real coefficients X, Y, X' and Y', which we thus come to put in place of the original ones of the bicomplex number, are nothing but the coefficients of the complex points X + iY, X' + iY lying respectively on the two proto-strings that contain the corresponding bicomplex point […]"*[35] [Segre 1892 p. 461], on which we will not dwell.

As already observed, moreover, Segre inserted his work in the more general context of the study of hypercomplex numbers systems by Weierstrass [Weierstrass 1884] and Dedekind [Dedekind 1885] and noticed that: *"The possibility of the decomposition of bicomplex numbers undertaken by us is a particular case of an analogous decomposition that Weierstrass and Dedekind undertake for hypercomplex numbers […]"*[36] [Segre 1892 p. 461].

Subsequently, bearing in mind once more the formulation by Weierstrass, he applied the decomposition found to the determination of the roots of an algebraic equation of degree $m$ to a bicomplex unknown [Segre 1892 p. 462]. He observed, precisely, that taking into account the aforesaid decomposition, this equation can be written $\sum_{l=0}^{m}(a_l + a'_l)(z+z')=0$, and be reduced, since the product of two "nullifics" of different sets is equal to zero, to $\sum_{l=0}^{m}a_l z^l + \sum_{l=0}^{m}a'_l z'^l = 0$, and therefore to the pair $\sum_{l=0}^{m}a_l z^l = 0$, $\sum_{l=0}^{m}a'_l z'^l = 0$. From this he therefore deduced that (if all the coefficients are not equal to zero) the equation has infinite solutions if and only if all its coefficients are zero divisors ("*nullifics*") of the same type (of the same set). Otherwise, denoting as $m$ and $m'$ the degrees of the two polynomials in $z$ and in $z'$, he deduced that the roots are $mm'$ and that if the given equation has degree $m$ and the coefficient of the term of degree $m$ is not a zero divisor ("*nullific*"), the equation has $m^2$ roots. Finally, Segre concluded this part with the analogous geometric result: *"In this way the preceding results are found analytically, that a straight line intersects a plane curve or complex algebraic surface of order m into $m^2$ bicomplex points, that two complex*

---

[34] "*I nullifiei della 1ª schiera sono dunque ridotti a numeri complessi (a coefficienti reali X, Y) con le due unità g, k; e similmente i nullifici della 2ª schiera con le due unità g', k' […]*"
[35] "*i nuovi coefficienti reali X, Y, X', Y', che così veniamo a sostituire a quelli primitivi del numero bicomplesso, non sono altro che i coefficienti dei punti complessi X + iY, X' + iY' giacenti risp. sui due protofili che contengono il corrispondente punto bicomplesso […]*"
[36] "*La possibilità della scomposizione dei numeri bicomplessi da noi operata rientra come caso particolare di un'analoga scomposizione che il Weierstrass ed il Dedekind fanno per numeri a quante si vogliano unità […]*"

*algebraic plane curves of orders m* and *$m_1$ intersect into $m^2 * m_1^2$ bicomplex points, etc."* [37][Segre 1892 p. 463]. Segre concluded the work hypothesizing further extensions, to *tricomplex numbers, s – complex numbers* and varieties *s – hyperalgebraic* ones, motivating their necessary introduction (paraphrasing Hamilton[38]) as follows: *"These subsequent extensions of the elements (geometric and analytical) and of the varieties formed with them appear in this way natural, spontaneous, and at the same time useful, indeed necessary. For ancient studies on algebraic entities real elements were sufficient: but for deeper studies it was necessary to introduce complex elements (ordinary). In the concept of these new elements there was certainly some discretion, and real elements could also have been generalized in other ways. The choice made was however the most appropriate for the algebraic field (in that for it an algebraic equation came to always have as many roots as is the degree, etc.). Once this choice is made the introduction of hyperalgebraic varieties, bicomplex elements, and so forth, is not an artifice anymore: it is, as we have said, a necessity."[39]* [Segre 1892 p. 465].

Probably the biggest lacuna in Segre's work was the absence of any reference to the British school, apart from Hamilton. This particularly refers to B. Peirce and W. Clifford, who are only marginally mentioned. Examination of their works would perhaps have made him go deeper into the difference between Hamilton's biquaternions and those of Clifford and notice the fact that his decomposition of bicomplex numbers in linear combination of idempotents is a particular case of that already pointed out in Peirce's 1870 work, published posthumously in 1881[40].

It is to be stressed that Segre's work was not immediately followed up, even among his direct alumni. It is to be noticed, however, that in 1898 an ex-alumnus of Segre's, Angelo Ramorino, devoted the last part of a historical work of his precisely to the study of bicomplex numbers.[41] We should also notice a late use of bicomplex numbers (in 1941) by one of Segre's greatest student, Gino Fano (1871 – 1952), with reference to an elementary geometry problem.[42]

---

[37] *"Pel tal modo si ritrovano analiticamente i risultati precedenti, che una retta taglia una curva piana o superficie algebrica d'ordine m complessa in $m^2$ punti bicomplessi, che due curve piane algebriche complesse di ordini m e $m_1$ si tagliano in $m^2 * m_1^2$ punti bicomplessi, ecc."*

[38] Segre refers to Hamilton before underlining the necessity of these extensions: *"To this, and particularly to the theory of bicomplex numbers (tricomplex ones, etc.) one can almost entirely apply the considerations that Hamilton made on biquaternions (we have already noticed the link between them and those numbers)"* [Segre 1892 p. 465]. (*"A ciò, e in particolare alla teoria dei numeri bicomplessi (tricomplessi ecc.) si possono applicare quasi completamente le considerazioni che Hamilton faceva sui biquaternioni (dei quali abbiam già rilevato il legame con quei numeri)"*) See the preceding section for Hamilton's considerations.

[39] *"Queste successive estensioni degli elementi (geometrici ed analitici) e delle varietà formate con essi appaiono per tal modo naturali, spontanee, ed in pari tempo utili, anzi necessarie. Per gli antichi studi degli enti algebrici bastavano gli elementi reali: ma per studi più profondi fu necessaria l'introduzione degli elementi complessi (ordinari). Nel concetto di questi nuovi elementi vi era certo dell'arbitrio, e gli elementi reali si sarebbero potuti generalizzare anche in altri modi. La scelta fatta fu però la più opportuna pel campo algebrico (in quanto per essa un'equazione algebrica veniva ad avere sempre tante radici quanto è il grado, ecc.). Una volta che tal scelta si è fatta l'introduzione delle varietà iperalgebriche, degli elementi bicomplessi, e così via, non è più un artificio: è, come abbiam detto, una necessità"*

[40] [Peirce 1881].

However, within Italian geometry (and above all analysis) the most important development to be noticed, in the 1930s and in the framework of the school of Gaetano Scorza (1876 – 1939), is the study of the analytical functions of bicomplex variables (and hypercomplex ones in general). Without going further into the matter, in this connection mention must be made of the works of Giuseppe Scorza Dragoni (1908 - 1996) and Nicolò Spampinato (1892-1971).[43]. This theme has had developments in fairly recent times, also in the framework of the study of the analytical functions of two complex variables.[44]

## 7. Conclusions

The analysis of the origins of algebra of "*bicomplex numbers*" confirms that in the middle of the nineteenth century in Great Britain there was great fervour around studies concerning new algebraic structures and around the meaning of symbolic algebra. These studies were then channelled into the general line of study of algebras that centred on the American school (A. Albert, L. Dickson, S. Epsteen, E. Moore, O. Veblen, J.M. Wedderburn) and the German one (R. Brauer, R. Dedekind, G. Frobenius, H. Hasse, E. Noether, E. Study). In comparison to this natural development, an exception in a sense is Segre, who introduced this subject in his studies on complex geometry, then studying their algebraic properties. In this connection, Segre believed that analysis could not only furnish some tools for geometry but that the opposite could also happen: "*hyperalgebraic entities, as I call them, were found in the analysts and not in the geometers, but vice versa it can happen that with the researches started by me and that I hope will be continued by other geometers, geometry will come in turn to be in the front line and furnish some new results that are useful for analysis. Do you not think so? In-depth synthetic study of a entity often leads to results which the analyst had not reached and that he will only reach after they are thus obtained by the geometer.*"[45] [Segre to Klein 1890 in [Luciano, Roero 2012]].

---

[41]On Angelo Ramorino (1869 -?) see the biography in http://www.peano2008.unito.it/scuola/ramorino.pdf. He had been an assistant of D'Ovidio until 1896 and in 1898 became an assistant of Peano. This seems to us one of the many connections between two opposite apparently schools like those of Segre and of Peano. The work mentioned is [Ramorino 1898].
[42][Fano 1941].
[43][Scorza Dragoni 1934]; [Spampinato 1935]; [Spampinato 1936].
[44]For a modern treatment of the matter from the algebraic point of view see, for instance: [Rochon, Shapiro 2004]; as regards aspects concerning complex analysis see [Luna-Elizarraras, Shapiro, Struppa, Vajiac 2012]. A general review of the state of the art on hypercomplex algebras can be found in [Olariu 2002].
[45] "*gli enti iperalgebrici com'io li chiamo, si trovavano negl'analisti e non nei geometri, viceversa può accadere che colle ricerche da me avviate e che spero vengano proseguite da altri geometri, la geometria venga a trovarsi a sua*

Specifically, the case of Cockle's "*Tessarines*", rediscovered and then called "*bicomplex numbers*" by Segre and within of the algebras of hypercomplex numbers, highlights the fact that a structure introduced in the context of strictly algebraic studies was rediscovered within complex geometry and found application in the study of the analytical functions of complex variables.

We also notice that some results of Segre's on the algebra of bicomplex numbers like isomorphism to the direct sum of two copies of the field of complex numbers and the fact that a polynomial equation of degree $n$ has $n^2$ roots (if they are not infinite) are inaccurately attributed to recent results, (1990s).[46]

Finally, further supporting what has been highlighted, we notice that the study of bicomplex numbers is today largely used in problems relating to generalizations of Mandelbrot sets and, closely connected to this, in the study of dynamic systems and, in physics, in quantum theory of fields and Solitons. But in relation to all of this we can only give some bibliographical references,[47] limiting ourselves to ascertaining the confirmation of what was foreseen by Segre on the "natural need" to introduce them[48] and in general maintained by Ian Stewart, according to whom the history of mathematics had repeatedly shown that discarding a fine and profound theory just because it had no immediate applications was a very bad move [Stewart 2008].

---

*volta in prima linea e fornire dei risultati nuovi ed utili all'analisi. Non le pare? Lo studio sintetico approfondito di un ente conduce spesso a risultati cui l'analista non era giunto e giungerà solo dopo che essi sono così ottenuti dal geometra*"
[46] See https://en.wikipedia.org/wiki/Bicomplex_number
[47] See for example [Rochon 2000]; [Matteau, Rochon 2015].
[48] In this connection see the preceding section.